\documentclass{article}
\usepackage[T1]{fontenc}

\usepackage{fullpage}
\newsavebox{\TableThree}
\usepackage{amsmath,amsthm,amssymb,xurl}
\usepackage{mathtools,pseudocode}
\usepackage{xurl}
\usepackage{enumitem} 
\usepackage{xspace} 
\usepackage{xifthen}
\usepackage[hidelinks]{hyperref}
\usepackage{tikz,authblk}
\usepackage{textcomp} 
\usepackage{fancybox}
\usepackage{float}
\usepackage{longtable}
\usepackage{threeparttablex}
\usepackage{authblk}
\usepackage{xcolor}
\usepackage{soul} 
\usepackage{multicol}
 \usepackage{array}
 \usepackage{ragged2e}
 \usepackage{caption}
\DeclareDocumentCommand{\F}{o}
{%
  \IfValueTF{#1}
    {\ensuremath{\textsc{f}(#1)}\xspace}
    {\ensuremath{\textsc{f}}\xspace}%
}

\newcommand{\zed}{\ensuremath{\mathbb{Z}}}
\newcommand{\Z}{\mathbb{Z}}

\newcommand{\T}{\ensuremath{\textsc{T}}}
\newcommand{\la}{\lambda}
\newcommand{\sm}{\setminus}
\newcommand{\Aut}{\ensuremath{\mathsf{Aut}}}
\newcommand{\can}{\ensuremath{\mathsf{canon}}}
\newcommand{\End}{\ensuremath{\mathsf{End}}}
\newcommand{\GL}{\ensuremath{\mathsf{GL}}}

\newcommand{\lexleq}{\ensuremath{\preceq}}
\newcommand{\lexless}{\ensuremath{\prec}}

\newcolumntype{R}[1]{>{\RaggedRight}p{#1}}

\theoremstyle{plain}
\newtheorem{theorem}{Theorem}[section]
\newtheorem{proposition}[theorem]{Proposition}
\newtheorem{lemma}[theorem]{Lemma}

\theoremstyle{definition}

\newtheorem{definition}[theorem]{Definition}

\newtheorem{example}[theorem]{Example}

{\refstepcounter{algorithm}\label{#1.alg}%
\stepcounter{Nalg}%
\begin{Sbox}\begin{minipage}{\boxwidth}%
{\bfseries Algorithm \thealgorithm:\tab{1}}\af{#1}($#2$)\\[2ex]%
\noindent\begin{math}\begin{array}{@{\tab{1}}l@{}}}%
{\end{array}\end{math}\end{minipage}%
\end{Sbox}\bigskip\noindent\fbox{\TheSbox}
\bigskip}

\floatstyle{plain}
\newfloat{algo}{htb}{.algo}

\allowdisplaybreaks

\usepackage{currfile}
\title{A computational approach to $(\lambda,\mu)$-fold factorizations of groups}

\author[1]{Donald L. Kreher}
\affil[1]{Department of Mathematical Sciences,
Michigan Technological University,
Houghton, MI 49931, U.S.A.}
\author[2]{Shuxing Li\thanks{Research supported by the U.S. National Science Foundation Grant DMS-2452236 and the University of Delaware Research Foundation Strategic Initiative (UDRF-SI) program.}}
\affil[2]{Department of Mathematical Sciences,
University of Delaware, Newark, DE 19716, U.S.A.}
\author[3]{Douglas R. Stinson}
\affil[3]{David R.\ Cheriton School of Computer Science,
University of Waterloo, Waterloo ON, N2L 3G1, Canada}

\date{}

\begin{document}
\maketitle

 \begin{abstract}
Let $G$ be a finite (multiplicative) group and let $S,T\subseteq G$. The subset pair $(S,T)$ is a \emph{$(\la,\mu)$-fold factorization} of $G$ if every nonidentity element of $G$ has exactly $\la$ representations as a product $xy$ with $x\in S$ and $y\in T$, and the identity element has exactly $\mu$ such representations. $(\la,\mu)$-fold factorizations provide a unified framework encompassing various combinatorial structures. Building on our prior theoretical study, we develop algorithms in this paper to find all inequivalent $(\la,\mu)$-fold factorizations for all abelian groups of order at most $34$. \end{abstract}
 
\section{Introduction}
\label{sec-introduction}

This paper is a sequel to our prior paper \cite{KLS-part1} establishing the theory of $(\la,\mu)$-fold factorizations of finite groups with $\mu \notin \{0, \la\}$. Here, we present a computational approach to investigating $(\la,\mu)$-fold factorizations of finite groups. By combining and extending algorithmic techniques described in 
\cite{KLS,KMS}, we obtain new examples and new nonexistence results for $(\la,\mu)$-fold factorizations through extensive computer searches.  This allows us to give a complete  classification of all inequivalent $(\la,\mu)$-fold factorizations of finite abelian groups of order at most $34$ (note that equivalence of factorizations is defined in Section \ref{sec-equivalence}).

\begin{definition}[{\cite[Definition 1.1]{KLS-part1}}]
\label{def-factorization}
Let $G$ be a finite (multiplicative) group and let $S,T\subseteq G$. The subset pair $(S,T)$ is a \emph{$(\la,\mu)$-fold factorization} of $G$ if every nonidentity element of $G$ has exactly $\la$ representations as a product $xy$ with $x\in S$ and $y\in T$, and the identity element has exactly $\mu$ such representations. The subsets $S$ and $T$ are the \emph{factors} in the factorization. If $|S|=s$ and $|T|=t$, the factorization has \emph{size $(s,t)$} and is denoted $\F[G,s,t,\la,\mu]$. 
\end{definition}

Various examples of factorizations can be found in Table \ref{Exists} in Section \ref{sec-computation}. Also, note that, unless an additive presentation is specified, a group $G$ is written multiplicatively with identity element $e$.

A subset of $G$ is identified with the sum of its elements in the group ring $\Z[G]$. Thus $(S,T)$ is a $(\la,\mu)$-fold factorization if and only if
\begin{equation}
\label{eqn-factorization}
ST=\la(G-e)+\mu e.
\end{equation}
We will denote $n=|G|$ throughout the paper. Therefore, the parameters of a factorization satisfy the following equation:
\begin{equation}
\label{eqn-counting}
st=\la(n-1)+\mu.
\end{equation}

Classical factorizations of groups, corresponding to $(1,1)$-fold factorizations, have been thoroughly studied; see the monographs \cite{Szabo04,SS}. \emph{Near-factorizations} of groups are $(1,0)$-fold factorizations, whose investigation dates back to de Bruijn \cite{DB56}. These factorizations have been further studied in \cite{BHS,CGHK,KMS,KPS2,Pech03,Pech04}. For $\la>1$, \emph{$\la$-fold near-factorizations} are $(\la,0)$-fold factorizations, which are proposed in \cite{KLS} and further studied in \cite{LM}. $(1,2)$-fold factorizations of groups have been studied in connection with thin Lehman matrices \cite{Pace,SS13,SS14,Shi}. Please refer to \cite[Section 1]{KLS-part1} for a more detailed account of prior research related to $(\la,\mu)$-fold factorization of finite groups.

\subsection{Organization of the paper}

In Sections \ref{sec-parameters}--\ref{sec-equivalence}, we review some basic theory and definitions from previous papers concerning group factorizations. In Section \ref{sec-automorphisms}, we provide a brief summary of the structure of  automorphism groups of abelian groups. In Section \ref{sec-inequivalent}, we describe the algorithms we use to generate inequivalent factorizations. We are mainly using the approaches described in \cite{KLS,KMS}. However, one useful optimization we employ in this paper is to use orderly algorithms to generate canonical factorizations (see Section \ref{canon.sec}). In Section \ref{sec-computation}, we report the results of our computer searches.

\subsection{Complements and basic parameter sets}
\label{sec-parameters}

Let $(S,T)$ be a $(\la,\mu)$-fold factorization of a group $G$. By taking complements of $S$ and/or $T$, we obtain a collection of four interrelated factorizations of $G$.

\begin{lemma}{\rm{(\cite[Lemma 3.2]{KLS-part1})}}
\label{lem-complement}
{The following statements 
are equivalent:}
\begin{enumerate}
\item $(S,T)$ is an $\F[G,s,t,\la,\mu]$.
\item $(G \sm S,T)$ is an $\F[G,n-s,t,t-\la,t-\mu]$.
\item $(S,G \sm T)$ is an $\F[G,s,n-t,s-\la,s-\mu]$.
\item $(G \sm S,G \sm T)$ is an $\F[G,n-s,n-t,n-s-t+\la,n-s-t+\mu]$.
\end{enumerate}
\end{lemma}
Together with interchange of the two factors, we can assume that $1\le s\le t\le n/2$. Following \cite[Section 3.1]{KLS-part1}, a \emph{basic parameter set} is an
integer tuple $(s,t,\la,\mu)$ satisfying
\begin{equation}
\label{eqn-basic}
\begin{split}
&st=\la(n-1)+\mu,\\
&1\le s\le t\le \frac{n}{2},\\
&1\le\mu\le s-1,\\
&\la\ne\mu.
\end{split}
\end{equation}
A factorization with a basic parameter set is called \emph{basic}. These conditions imply $s\ge 2$ and $1\le\la\le s-1$ \cite[Remark 3.3]{KLS-part1}. The cases $\mu=0$ and $\mu=s$ are related
by complementation, where $\mu=0$ corresponds to $\la$-fold near-factorizations. As observed in \cite[Remark 3.8]{KLS-part1}, the case $\la=\mu$ is very  different from $(\la,\mu)$-fold factorizations with $\la \ne \mu$. From now on, we always consider basic $(\la,\mu)$-fold factorizations.

The classification problem considered here is to determine, for every abelian group $G$ of order at most $34$ and every basic parameter set $(s,t,\la,\mu)$, all $\F[G,s,t,\la,\mu]$ up to the equivalence defined
in Section~\ref{sec-equivalence}.

\subsection{Group matrices and uniqueness of mates}
\label{sec-matrices}

\begin{definition}[{\cite[Definition 3.4]{KLS-part1}}]
\label{def-groupmatrix}
An $n\times n$ real-valued matrix $M$ whose rows and columns are indexed by $G$ is a
\emph{group matrix} based on $G$ if
\[
M_{gx,gy}=M_{x,y}\qquad\text{for all }g,x,y\in G.
\]
Equivalently, the $(x,y)$ entry of $M$ depends only on the quotient $x^{-1}y$.
\end{definition}

For $A\subseteq G$, the group matrix associated with $A$ is $M(A)$, where 
\begin{equation}
\label{eqn-mat}
\bigl(M(A)\bigr)_{g,h}=
\begin{cases}
1 & \text{if }g^{-1}h\in A,\\
0 & \text{otherwise}.
\end{cases}
\end{equation}
With a fixed ordering of group elements, the subset $A$ is uniquely determined by $M(A)$.
The row and column sums of $M(A)$ are $|A|$. We use $I$ and $J$ for the
$n\times n$ identity matrix and all-ones matrix, respectively.

\begin{theorem}{\rm (\cite[Theorem 3.5]{KLS-part1})}
\label{thm-matrixrep}
Suppose $S,T \subseteq G$. Then $(S,T)$ is a $(\la,\mu)$-fold factorization of $G$ if and only if 
\begin{equation}
\label{eqn-matrix}
M(S)M(T)=\la (J-I) + \mu I.
\end{equation}
\end{theorem}

For a basic parameter set, the right-hand side of (\ref{eqn-matrix}) has determinant 
\[(\mu-\la)^{n-1}\big(\la n+\mu-\la \big)=(\mu-\la)^{n-1}st\ne0.\] Hence $M(S)$ and $M(T)$ are nonsingular, because $\lambda \neq \mu$. 

This matrix perspective leads to the following results. 

\begin{theorem}{\rm (\cite[Theorem 3.7]{KLS-part1})}
Let $G$ be a finite group. Let $\la$ be a positive integer and let $\mu$ be a nonnegative integer with $\la \ne \mu$.
\label{thm-unique}
\begin{enumerate}
\item Suppose $(S,T)$ and $(S',T)$ are both $(\la,\mu)$-fold factorizations of $G$. Then $S = S'$.
\item Suppose $(S,T)$ and $(S,T')$ are both $(\la,\mu)$-fold factorizations of $G$. Then $T = T'$.
\end{enumerate}
\end{theorem}
Consequently, for fixed $G,\la$ and $\mu$ with $\la \ne \mu$, either factor  uniquely determines the other. Therefore, we say one factor is the \emph{mate} of the other. 

\smallskip

Another, stronger, uniqueness property has been proven in \cite{KLS-part1}.

\begin{theorem}{\rm (\cite[Theorem 3.10]{KLS-part1})}
\label{thm-unique2}
Suppose $(S,T)$ is a $(\la,\mu)$-fold factorization of $G$ and $(S',T)$ is a $(\la',\mu')$-fold factorization of $G$, where $S \neq S'$, $\la \neq \mu$ and $\la' \neq \mu'$. Suppose $1 \leq s,s' \leq |G|-1$, where $s = |S|$ and $s' = |S'|$. Then $S' = G \sm S$, $\la' = t - \la$ and $\mu' = t - \mu$, where $|T| = t$.
\end{theorem}

\subsection{Equivalence of factorizations and canonical factorizations}
\label{sec-equivalence}

For the purpose of classification, we need to define equivalence and canonical forms of subsets and factorizations.  First, we fix a total ordering $e = g_1 < g_2 < \cdots <g_n$ of the elements of $G$. {This ordering will be used throughout the rest of the paper.} 
Elements in a subset  are written in increasing order. 

Subsets of the same size are compared in lexicographic order. {Suppose that 
$S = \{ s_i : 1 \leq i \leq m\}$ and $T = \{ t_i : 1 \leq i \leq m\}$, where
$s_1 < \cdots < s_m$ and $t_1 < \cdots < t_m.$ We say that $S \lexless T$ if there exists an index $j$, where $0 \leq j\leq m-1$,  such that
$s_i = t_i$ for $1 \leq i \leq j$ and $s_{j+1} < t_{j+1}$. Also, $S \lexleq T$  if $S \lexless T$ or $S = T$.
}

If $|S| < |T|$, then we use the convention that $S \lexless T$. 
Thus we have a total ordering on all the possible subsets of $G$.
We can also extend this ordering $\lexleq$ to ordered pairs of subsets. We say that
$(S,T) \lexleq (S',T')$ if $S \lexless S'$, or $S = S'$ and $T \lexleq T'$.

\smallskip

We now define equivalence of subsets. Let $\Aut(G)$ denote the automorphism group of $G$.

\begin{definition}
\label{subsest-equivalence}
Let $S$ and $S'$ be subsets of a finite abelian group $G$, where $|S| = |S'|$.
We say that $S$ and $S'$ are \emph{equivalent} if there exist
$f\in\Aut(G)$ and $g\in G$ such that
$
S'=f(S)g$.
\end{definition}

We will use the notation $\sim$ to denote equivalence of subsets. It is not hard to prove that 
$\sim$ is indeed an equivalence relation. The equivalence class containing a subset $S$ will be denoted
by $[S]$.

\begin{definition}
\label{def-canonical-subset}
For $A\subseteq G$, define
$$
\can(A)
=\min_{\substack{f\in\Aut(G)\\g\in G}} f(A)g.
$$
The subset $A$ is \emph{canonical} if $A=\can(A)$.
\end{definition}

Thus $A$ is canonical if and only if $A\lexleq f(A)g$ for every $f\in\Aut(G)$ and every $g\in G$.
Equivalently, $A$ is the smallest subset in the equivalence class $[A]$.

\smallskip

The equivalence notion $\sim$ extends to ordered pairs of subsets and, in particular, factorizations. But we need to distinguish between the two cases when the subsets have different sizes or the same sizes.

\begin{definition}{\rm{(\cite[Definition 3.15]{KLS-part1})}}
\label{def-equivalence}
Let $G$ be an abelian group and let $(S,T)$ and $(S',T')$ both be $\F[G,s,t,\la,\mu]$.
If $s<t$, the factorizations $(S,T)$ and $(S',T')$ are \emph{equivalent} if there exist
$f\in\Aut(G)$ and $g\in G$ such that
$$
S'=f(S)g \quad \text{and} \quad T'=f(T)g^{-1}.
$$
If $s=t$, the factorizations $(S,T)$ and $(S',T')$ are equivalent if there exist $f\in\Aut(G)$ and $g\in G$
such that either
$$
(S',T')=(f(S)g,f(T)g^{-1})
$$
or
$$
(S',T')=(f(T)g,f(S)g^{-1}).
$$
\end{definition}
Note that, when $|S| = |T|$, we also allow for the possibility of interchanging $S$ and $T$ in the above definition.

\smallskip

We again use the notation $\sim$ to denote equivalence of factorizations. 
The equivalence class containing a factorization $(S,T)$ will be denoted by
$[(S,T)]$. We define canonicity of factorizations as follows.

\begin{definition}
\label{def-canonical-factorization}
An $\F[G,s,t,\la,\mu]$, say $(S,T)$, is \emph{canonical} if one of the following conditions holds:
\begin{enumerate}
\item $s<t$ and $S=\can(S)$; or
\item $s=t$, $S=\can(S)$ and $S\lexleq\can(T)$.
\end{enumerate}
\end{definition}

We define $\can(S,T)$  to be the smallest
factorization in the equivalence class $[(S,T)]$. Thus $(S,T) = \can(S,T)$ precisely when $(S,T)$ is canonical.

\smallskip
{The following example illustrates the extra work that needs to be done to compute $\can(S,T)$ when $|S| = |T|$.}

\begin{example}
Let $G = \zed_{21}$; $G$ is an additive abelian group. 
The units in $\zed_{21}$ are denoted by
 \[\zed_{21}^* = \{ 1,2,4,5,8,10,11,13,16,17,19,20 \}.\] 
 It is not hard to see that 
\[\Aut(G) = \{ x \mapsto mx \bmod 21: m \in  \zed_{21}^* \}.\] 
The transformations we use to define equivalence have the form 
$x \mapsto mx + g \bmod 21$, where $m \in \zed_{21}^*$ and $g \in \zed_{21}$.

Define
$S=\{0,1,2,3,5,7,17,19\}$ and $T=\{0,7,8,10,11,18,19,20\}$. 
It can be verified that $(S,T)$ is an $\F[G,8,8,3,4]$.

It is  straightforward to check that $S$ is canonical.
However, it turns out that $(S,T)$ is not canonical. This is because 
 $\can(T)=\{0,1,2,3,5,6,7,8\}$, obtained by applying the mapping $x \mapsto 2x+7 \bmod 21$ to $T$. 
 We see that $\can(T) \lexless \can(S) = S \lexless T$. 
 The canonical form of $(S,T)$ is $(T',S')$, where
 \[ T' = 2T+7 = \{0,1,2,3,5,6,7,8\}\]
 and
 \[S' = 2S-7 = \{3,6,7,10,14,16,18,20 \}.\tag*{$\blacksquare$}\]
\end{example}

The following result is straightforward. 

\begin{proposition}
\label{prop-canonical}
Suppose $G$ is an abelian group, $(s,t,\la,\mu)$ is a basic parameter set and $(S,T)$ is an $\F[G,s,t,\la,\mu]$.
Then $[(S,T)]$ contains exactly one canonical factorization.
\end{proposition}
\begin{proof}
Suppose first that $s<t$. Choose $f\in\Aut(G)$ and $g\in G$ such that $f(S)g=\can(S)$. The equivalent factorization 
$$
(f(S)g,f(T)g^{-1})=(\can(S),f(T)g^{-1})
$$ 
is canonical from Definition \ref{def-canonical-factorization}. 

Assume that there are two  canonical factorizations equivalent to $(S,T)$. These two factorizations must have the same first factor, namely, $\can(S)$. Hence, by Theorem \ref{thm-unique}, they have the same second factor. Therefore, there exists a unique canonical factorization within an equivalence class. 

Now we consider the case $s=t$. 
First, suppose that $\can(S) \lexleq \can(T)$.

Let $S' = \can(S) = f(S)g$ and define $T' =  f(T)g^{-1}$.
We will show that $(S',T')$ is canonical.
We have $(S',T') \sim (S,T)$ and $\can(S) = \can(S')$ because $S \sim S'$. 
Since $S' = \can(S)$, we have $\can(S') = S'$.

Hence, from Definition \ref{def-canonical-factorization}, in order to prove canonicity of $(S',T')$, it suffices to show that $S' \lexleq \can(T')$.
We have assumed that $\can(S) \lexleq \can(T)$. 
Also, $\can(T) = \can(T')$ because $T \sim T'$.
Therefore,
\[S' = \can(S') = \can(S) \lexleq  \can(T) = \can(T'),\] as desired.

Now assume that there are two canonical factorizations equivalent to $(S,T)$.
These two factorizations must have the same first factor, namely, $\can(S)$. Hence, by Theorem \ref{thm-unique}, they have the same second factor. Therefore, there exists a unique canonical factorization within an equivalence class. 

If $\can(T) \lexleq \can(S)$, then the proof is similar; we omit the details.
\end{proof}

\section{Automorphisms of finite abelian groups}
\label{sec-automorphisms}

We extensively employ automorphism groups to prune our factorization search algorithms in Section~\ref{sec-inequivalent}. This section provides a brief overview of the computation of the automorphism group of an abelian group. For additional details, please refer to \cite{dummitfoote,Hillar-Rhea,Mader}. Throughout this section, the group operation is written additively.

According to \cite[Lemma 2.1]{Hillar-Rhea}, if $H$ and $K$ are finite groups of relatively prime order, then  
$$
\Aut(H \times K) \cong \Aut(H) \times \Aut(K).  
$$ 
Consequently, because each finite abelian group decomposes as the direct product of its Sylow $p$-groups, it suffices to consider finite abelian $p$-groups of the form  
$$
G_p = \Z_{p^{m_1} } \times \cdots \times \Z_{p^{m_r}},  
$$
where $p$ is prime and $m_1 \geq \cdots \geq m_r \geq 1$. 

First, we present some definitions and notation.
\begin{itemize}
\item Let $\End(G_p)$ be the group of homomorphisms from $G_p$ to itself, with the group operation given by $(A + B)(g)= A(g)+B(g)$ for $A, B \in \End(G_p)$ and $g \in G_p$.
\item Let $R_p=\{ (a_{ij}) \in \Z^{r \times r}:  p^{\max\{m_i-m_j, 0\}} \mid a_{ij}, \forall 1 \le i,j \le r\}$.
Thus $R_p$ is a set of $r$ by $r$ integer matrices whose entries satisfy certain divisibility conditions.
\item For an integer vector $\mathbf{g} = (g_1,g_2,\ldots,g_r)^T$ of length $r$, let $\pi (\mathbf{g})$ denote the vector obtained by reducing the $i$th component of  $\mathbf{g}$ modulo $p^{m_i}$, for $1 \leq i \leq r$.
\item Define the map $\Psi : R_p \to \End(G_p)$ by
 $$
 \Psi(A)((g_1,g_2,\ldots,g_r)^T)=\pi(A(g_1,g_2,\ldots,g_r)^T).
 $$ Thus, a mapping $\Psi (A)$ performs a matrix-vector multiplication followed by the modular reductions specified by $\pi$. We observe that $\Psi (A)$ is an endomorphism precisely because of the divisibility conditions imposed on the entries of a matrix 
 $A \in R_p$. Basically, we are making use of the fact that
 $\zed_{p^{m_j}}$ is embedded in $\zed_{p^{m_i}}$ when $m_i > m_j$ via the mapping
 $x \mapsto  p^{m_i-m_j}x$.
\item $\GL(r,\zed_p)$ is the general linear group of degree $r$ over the field $\zed_p$. That is, $\GL(r,\zed_p)$ consists of all the $r$ by $r$ invertible matrices with entries from $\zed_p$.
\end{itemize}

\begin{lemma}
{\rm (\cite[Theorem 3.3, Lemma 3.4]{Hillar-Rhea})}
The map $\Psi : R_p \to \End(G_p)$ 
is a surjective homomorphism with kernel
$$
\ker \Psi=\{ (a_{ij}) \in R_p:  p^{m_i} \mid a_{ij} \text{ for all } 1 \le i,j \le r\}.
$$
Therefore, $\End(G_p) \cong R_p/\ker \Psi$.
\end{lemma}

\begin{lemma}
{\rm (\cite[Theorem 3.6]{Hillar-Rhea})}
Let $A \in R_p$. The  endomorphism $\Psi(A)$ is an automorphism of $G_p$ if and only if $A \bmod p \in \GL(r,\zed_p)$. Consequently, 
$$
\Aut(G_p)=\{ \Psi(A): A \in R_p,  A \bmod p \in \GL(r,\zed_p) \}.
$$ 
\end{lemma}

For $A \in R_p$, the matrix $A \bmod p$ has a nice block structure. Consider the multiset of exponents,
$\{ m_1, \dots , m_r \}$. Denote the multiplicities of the elements of this multiset by
$n_1, \dots , n_s$, where $\sum n_i = r$. Denote the (square) blocks of sizes $n_1, \dots , n_s$ on the diagonal of $A \bmod p$ by $A_1, \dots , A_s$. All of the elements of $A \bmod p$ above these blocks are $0$s. Hence, $\det (A \bmod p) \not\equiv 0 \bmod p$ if and only if $\det (A_i) \not\equiv 0 \bmod p$ for $1 \leq i \leq s$. This happens if and only if $A_i \in \GL(n_i,\zed_p)$ for $1 \leq i \leq s$.

Due to this nice block structure, it is possible to give a formula for $|\Aut(G_p)| $.

\begin{lemma}
{\rm (\cite[Theorem 4.1]{Hillar-Rhea})}
\label{autsize.lem}
For $1 \le k \le r$, define 
$$
d_k=\min \{ 1 \le \ell \le k : m_\ell=m_k \}, \quad c_k=\max \{ k \le \ell \le r: m_\ell=m_k \}.
$$
Then
\begin{equation}
\label{autsize.eq}
|\Aut(G_p)| = \prod_{k=1}^r (p^{r+1-d_k}-p^{r-k}) \prod_{j=1}^r (p^{m_j})^{d_j-1} \prod_{i=1}^r  (p^{m_i-1})^{c_i}.
\end{equation}
\end{lemma}

We have the following example illustrating the above theorems.

\begin{example}
Let $G=\Z_8 \times \Z_2 \times \Z_2$. Thus, following the notation in Lemma \ref{autsize.lem}, we have $p=2$, $r=3$, $m_1=3$ and $m_2=m_3=1$.  Then we have  $(d_1,d_2,d_3)=(1,2,2)$ and $(c_1,c_2,c_3)=(1,3,3)$. Consequently,
\begin{align*}
|\Aut(G)|&=\prod_{k=1}^3 (2^{4-d_k}-2^{3-k}) \prod_{j=1}^3 (2^{m_j})^{d_j-1} \prod_{i=1}^3  (2^{m_i-1})^{c_i} \\
             &=(2^3-2^2)(2^2-2)(2^2-1) \cdot (2^3)^{1-1}(2^1)^{2-1} (2^1)^{2-1}  \cdot (2^2)^{1} (2^0)^{3} (2^0)^{3}=384.
\end{align*}
Note that 
$$
R_2=\left\{ \begin{pmatrix}  
              b_{11} & 4b_{12} & 4b_{13} \\
              b_{21} & b_{22} & b_{23} \\
              b_{31} & b_{32} & b_{33}
             \end{pmatrix}: b_{ij} \in \Z  \right\}.
$$
Let $\Psi$ be the homomorphism from $R_2$ to $\End(G)$. 
Then
$$
\ker \Psi=\left\{ \begin{pmatrix}  
              8b_{11} & 8b_{12} & 8b_{13} \\
              2b_{21} & 2b_{22} & 2b_{23} \\
              2b_{31} & 2b_{32} & 2b_{33}
             \end{pmatrix}: b_{ij} \in \Z  \right\}.
$$
Therefore, a complete set of coset representatives of $\ker \Psi$ in $R_2$ is
$$
\End(G) \cong \left\{\Psi(A):  A = \begin{pmatrix}  
              b_{11} & 4b_{12} & 4b_{13} \\
              b_{21} & b_{22} & b_{23} \\
              b_{31} & b_{32} & b_{33}
             \end{pmatrix}: b_{11} \in \{0,1,\dots, 7\}, b_{12}, b_{13}, b_{21}, b_{22}, b_{23}, b_{31}, b_{32}, b_{33} \in \{0,1\} \right\}.
$$

Specifically, each coset representative $A$
gives rise to a homomorphism $\Psi(A) \in \End(G)$ such that, for $x=(x_1,x_2,x_3)^T \in \Z_8 \times \Z_2 \times \Z_2$, we have 
$$
\Psi(A)((x_1,x_2,x_3)^T)=\pi(A(x_1,x_2,x_3)^T)=\begin{pmatrix}  
              b_{11}x_1+4b_{12}x_2+4b_{13}x_3 \bmod 8 \\
              b_{21}x_1+b_{22}x_2+b_{23}x_3 \bmod 2\\
              b_{31}x_1+b_{32}x_2+b_{33}x_3 \bmod 2
     \end{pmatrix},
$$
where $A(x_1,x_2,x_3)^T$ is regarded as a matrix-vector product over $\Z$.
Moreover, $\Psi(A) \in \Aut(G)$ if and only if $A \bmod 2 \in \GL(3,\zed_2)$. This happens if and only if
\begin{enumerate}
\item $b_{11}$ is odd, and
\item $\det \left(
\begin{array}{cc}
b_{22} & b_{23} \\
 b_{32} & b_{33}
\end{array}
\right)$ is odd.
\end{enumerate}
There are four ways to choose $b_{11} \in \{0,1,\ldots,7\}$, two ways to choose each of $b_{12}$, $b_{13}$, $b_{21}$ and  $b_{31}$ (which are elements of $\{0,1\}$), and six ways to construct the submatrix $\left( 
\begin{array}{cc}
b_{22} & b_{23} \\
 b_{32} & b_{33}
\end{array}
\right)$, where $b_{22}, b_{23},
 b_{32}$  and $b_{33}$ are elements of $\zed_2$. This yields a total of $4 \times 2^4 \times 6 = 384$ matrices, agreeing with 
equation (\ref{autsize.eq}).\hfill$\blacksquare$
\end{example}

\section{Enumerating inequivalent factorizations}
\label{sec-inequivalent}

There are several aspects to be considered in the enumeration of all inequivalent factorizations $\F[G,s,t,\la,\mu]$, given an abelian group $G$ and a basic parameter set $(s,t,\la,\mu)$.
\begin{enumerate}
\item Canonicity testing and computing canonical forms
\item Finding mates
\item Generating canonical factorizations.
\end{enumerate}
We consider each of these three steps in turn in the following subsections.

Note that we are assuming here that the basic parameter set $(s,t,\la,\mu)$ is specified ahead of time. In Section \ref{Sec3.4}, we will consider a modified setting where only the group $G$ is fixed and we want to  test all basic parameter sets.

\subsection{Canonicity testing}
\label{canonicaltest.sec}

In Algorithm \ref{CCS}, we present pseudocode 
that computes $\can(S)$, given a subset $S \subseteq G$ as input.
We assume that $\Aut(G)$ has already been computed, e.g.,  using techniques described in Section \ref{sec-automorphisms}. In Algorithm \ref{CCS}, we use the ordering $\lexleq$ that was defined in Section \ref{sec-equivalence}. Also we denote the identity {automorphism} in $\Aut(G)$ by $\mathsf{id}$.
Algorithm \ref{CCS} returns a triple $(S',f',g')$, where $S' = \can(S) = f'(S)g'$. 

\begin{algo}
\begin{pseudocode}[shadowbox]{ComputeCanonicalSubset}{S}
\label{CCS}
\GLOBAL{G,\Aut(G)}\\
S' \GETS S\\
f' \GETS \mathsf{id}\\
g'\GETS e\\
\FORALL g \in G \DO \BEGIN
\FORALL f \in \Aut(G) \DO \BEGIN
S'' \GETS f({S})g\\
\IF S'' \lexless S' \THEN  
\BEGIN 
S' \GETS S''\\
f' \GETS f\\
g' \GETS g\\
\END
\END\\
\END\\
 \RETURN{S',f',g'}
\end{pseudocode}
\end{algo}

In Algorithm \ref{CCF}, we present pseudocode 
that computes $\can(S,T)$, given a factorization $(S,T)$ as input.
In this algorithm, we assume that $(S,T)$ is a factorization, $S$ is canonical and $|S| \leq |T|$.
The correctness of Algorithm \ref{CCF} follows from Proposition \ref{prop-canonical}.

\begin{algo}
\begin{pseudocode}[shadowbox]{ComputeCanonicalFactorization}{S,T}
\label{CCF}
\COMMENT{We assume that $(S,T)$ is a factorization, $|S| \leq |T|$ and $S$ is canonical}\\
\EXTERNAL{\mathsf{ComputeCanonicalSubset}}\\
\GLOBAL{G,\Aut(G)}\\
\IF |S| < |T| \THEN
\RETURN{S,T}
\ELSE
\BEGIN 
\COMMENT{here we have $|S| = |T|$}\\
(T',f',g') \GETS \mathsf{ComputeCanonicalSubset}(T)\\
{\COMMENT{$T' = f'(T)g'$}}\\
\IF T' \lexless S \THEN \BEGIN
S' \GETS f'(S){(g')^{-1}}\\
\RETURN {T',S'} \END
\ELSE 
\RETURN {S,T}
\END 
\end{pseudocode}
\end{algo}

\bigskip

\subsection{Finding mates}
\label{mates.sec}

Given $S$, we want to find a mate $T$. There are two (equivalent) approaches:
\begin{description}
\item[method 1] \mbox{\quad} \\ Compute the matrix $M(S)^{-1}$ (if it exists) and then use equation (\ref{eqn-matrix}) to find $M(T)$. If $M(T)$ is a 0-1 matrix, then we have a mate $T$ for $S$.
\item[method 2]  \mbox{\quad} \\ Solve an associated system of linear equations to see if there is a 0-1 solution. Specifically, solve the linear system $M(S)(x_1, \dots , x_n)^T = (\mu,\lambda, \dots , \lambda)^T$, {where $M(S)$ is a 0-1 matrix and $\lambda, \mu$ are positive integers}. If there is a 0-1 solution $(x_1, \dots , x_n)$, then we can easily retrieve the mate $T$.
\end{description}

Both of these methods are described in detail in previous papers (in particular, see \cite{KLS,KMS}). The paper \cite{KMS} discussed the case of {$(1,0)$-fold factorizations} and \cite{KLS} addressed finding mates that yield {$(\lambda,0)$-fold factorizations}. The methods generalize in a straightforward manner to {$(\lambda,\mu)$-fold factorizations} when $\lambda \neq \mu$.

We have found in practice that using ``Dancing Links'' \cite{KP,Knuth} is an efficient way to solve the system of linear equations mentioned above as method 2. We note that the Dancing Links approach requires that $\mu \ne 0$, {so additional preprocessing was required when dealing with $(1,0)$-fold factorizations \cite{KMS} and $(\lambda,0)$-fold factorizations \cite{KLS}, where $\mu=0$.} However, that is not an issue here, because we are assuming $\mu \neq 0$.

\subsection{Generating canonical factorizations}
\label{canon.sec}

For a specified abelian group $G$ and a basic parameter set $(s,t,\lambda,\mu)$, we want to generate all 
inequivalent $\F[G,s,t,\la,\mu]$, where equivalence is defined in Definition \ref{def-equivalence}. 
In view of Proposition \ref{prop-canonical},
it suffices to generate all canonical $\F[G,s,t,\la,\mu]$. All the algorithms are based on the following general approach.

\begin{description}
\item[step 1]  Generate all the canonical $s$-subsets $S$.
\item[step 2]  For each such set $S$ generated in step 1, determine if $S$ has a mate $T$. If so, proceed to step 3.
\item[step 3]  Test $(S,T)$ to see if it is canonical. If so, then output $(S,T)$. 
\end{description}

We observe that, in step 3, we are guaranteed that $S$ is canonical. Therefore we can use the algorithm $\mathsf{ComputeCanonicalFactorization}$, which was defined in Algorithm \ref{CCF}, to test canonicity: the factorization $(S,T)$ is canonical if and only if $(S,T) = \mathsf{ComputeCanonicalFactorization}(S,T)$.

Steps 2 and 3 have been addressed already, in Sections \ref{mates.sec} and \ref{canonicaltest.sec}, respectively. There are various ways to implement step 1, which we discuss now.

\smallskip

In the \emph{basic method}, we
generate all the $s$-subsets $S$ that contain $e$ (the identity). This can be done using any algorithm to generate all the subsets of a given size of a given set. Specifically, we would compute $S = \{e \} \cup S'$ for all $(s-1)$-subsets $S' \subseteq G \setminus \{e\}$. Then, for every $S$, we test to see if it is canonical. If so, we proceed to  step 2.

However, it is much more efficient to use an \emph{orderly algorithm} (see Read \cite{Read}).
The defining feature of an orderly algorithm is that we test every $j$-subset for canonicity, for $1 \leq j \leq |S|$. (As we will see, if a subset $S$ is canonical, then the $j$ smallest elements of $S$ form a canonical subset, for all $j$ such that $1 \leq j \leq |S|$.) It is obviously more work to test partial ``structures'' for canonicity, but there is a huge savings, as the resulting pruning of the search tree is very effective.
More formally, the correctness of 
the orderly method is a consequence of the following lemma.

\begin{lemma}
\label{lem-canonical-prefix}
{Suppose $a_1<\cdots<a_k$ and $A=\{a_1, \dots , a_k\}$ is canonical.} Then
$\{a_1,\ldots,a_j\}$ is canonical for every $1\le j\le k$.
\end{lemma}
\begin{proof}
Suppose $B = \{a_1,\ldots,a_j\}$ is not canonical, for some $j$, $1\le j\le k$.
Let $B' = \can(B)$; then $B' \lexless B$. Suppose that $B' = f(B)g$ where $f \in \Aut(G)$ and $g \in G$.
Let $A' = f(A)g$ and let $C$ denote the $j$ smallest elements in $A'$.
Then $C \lexleq B'$, so $C \lexless B$. Hence, $A' \lexless A$, which contradicts the assumption that $A$ is canonical.
\end{proof}

In designing an orderly algorithm, there are two approaches, termed \emph{breadth-first} and \emph{depth-first}.
In the breadth-first method, we create a list $\mathcal{L}_j$ of all canonical $j$-subsets, for $j = 1, \dots , s$. Given the list $\mathcal{L}_j$, it is straightforward to compute $\mathcal{L}_{j+1}$. Once we have created $\mathcal{L}_s$, we proceed to check every $S  \in \mathcal{L}_s$ to see if it has a mate. 
On the other hand, in the depth-first approach, we use a recursive algorithm to generate the canonical $s$-subsets. Every time we construct a (canonical) $s$-subset $S$, we check to see if it has a mate.
We present more detailed descriptions of both approaches in the following subsections.

\subsubsection{Breadth-first approach}

In a given iteration of  the breadth-first algorithm, 
we have a list $\mathcal{L}_j$ of all canonical $j$-subsets of $G$, in increasing order. 
Then we generate all the canonical $(j+1)$-subsets of {$G$}. Initially, $j = 1$ and the list $\mathcal{L}_1$ consists of the single subset {$\{g_1= e\}$}.

For a given $j$, we iterate through the sets in the list $\mathcal{L}_j$. 
For each set $S \in \mathcal{L}_j$, we can add every possible element to $S$ 
(noting that we only need to add elements that are greater than the largest element of $S$). Each of the resulting sets $S'$ is tested for canonicity. If $S'$ is canonical, then it is appended to the new list $L_{j+1}$.
The final list $\mathcal{L}_s$ is the desired list of canonical $s$-subsets.
In Algorithm \ref{BFCS}, we present a pseudocode description of this process.

Note that Algorithm \ref{BFCS} returns all of the canonical subsets of size at most $s$ in $s$ separate lists; this will be useful in Section \ref{Sec3.4}. However, for the discussion in the current section, we only need $\mathcal{L}_s$.

\begin{algo}
\begin{pseudocode}[shadowbox]{BreadthFirstCanonicalSubsets}{s}
\label{BFCS}
\GLOBAL{{G=\{g_1=e,g_2,\ldots,g_n\}},\Aut(G)}\\
\mathcal{L}_1 \GETS [\{g_1\}]\\
 \FOR j \GETS 1 \TO s-1 \DO \BEGIN
\mathcal{L}_{j+1} \GETS [\: ] \quad \text{(an empty list)}\\
\FORALL S \in \mathcal{L}_{j} \DO \BEGIN
 \text{denote the maximum element of $S$ by $g_m$}\\
  \FOR i \GETS m+1 \TO n \DO \BEGIN
 S' \GETS S \cup \{g_i\}\\
 \IF S' {\text{ is a canonical subset}} \CTHEN \\
 \text{append $S'$ to $\mathcal{L}_{j+1}$}\\
\END
\END
\END \\
\RETURN{\mathcal{L}_{1}, \dots , \mathcal{L}_{s}}
\end{pseudocode}
\end{algo}

{The following example records a collection of canonical subsets of increasing size being generated sequentially.}

\begin{example}
Suppose $G = \zed_{12}$ and $s = 4$.
We order the group elements in the natural way: $0 < \cdots < 11$.
We have the following lists of canonical subsets of various sizes:
\[ \mathcal{L}_{1} = [\{0\}],\]

\[\mathcal{L}_{2} = [\{0,1\},\{0,2\},\{0,3\},\{0,4\},\{0,6\}],\]

\[\mathcal{L}_{3} = [\{0, 1, 2\}, \{0, 1, 3\}, \{0, 1, 4\}, \{0, 1, 5\}, \{0, 1, 6\}, \{0, 2, 4\}, \{0, 2, 6\}, \{0, 3, 6\}, \{0, 4, 8\}],\]
and
\begin{align*}
\mathcal{L}_{4} =[& \{0,1,2,3\},\{0,1,2,4\},\{0,1,2,5\},\{0,1,2,6\},\{0,1,2,7\}, \{0,1,3,4\},\{0,1,3,5\}, \\
                            & \{0,1,3,6\},\{0,1,3,7\},\{0,1,3,9\},\{0,1,3,10\},\{0,1,4,5\},\{0,1,4,7\},\{0,1,4,8\}, \\
                            & \{0,1,4,9\},\{0,1,5,6\},\{0,1,6,7\},\{0,2,4,6\},\{0,2,4,8\},\{0,2,6,8\},\{0,3,6,9\} ].
\tag*{$\blacksquare$}
\end{align*}
\end{example}

\bigskip

{The following example illustrates the reduction achieved by restricting our focus to canonical subsets.}

\begin{example}
The number of subsets in these lists is reduced considerably by requiring each list to consist only of canonical subsets. We present some data pertaining to $\zed_{30}$. We indicate the reduction, for subsets of sizes $\leq 8$.
\[
\begin{array}{c|c|c}
j & \text{number of $j$-subsets} & \text{number of canonical $j$-subsets}\\
\hline
1 & 30 & 1\\
2 & 435 & 7\\
3 & 4060 & 25\\
4 & 27405 & 149\\
5 & 142506 & 646\\
6 & 593775 & 2643\\
7 & 2035800 & 8708\\
8 & 5852925 & 24938
\end{array}
\]

\smallskip

It is clear that there is a significant reduction when restricting to canonical subsets.\hfill$\blacksquare$
\end{example}

\subsubsection{Depth-first approach}

The depth-first approach is based on a recursive backtracking algorithm, which we name $\mathsf{Extend}$.
We end up generating all the canonical $s$-subsets in lexicographical order, just as in the breadth-first approach.
However, in the depth-first algorithm, whenever we construct a canonical $s$-subset $S$, we proceed immediately to step 2 to see if $S$ has a mate. In the breadth-first approach, we place all the canonical $s$-subsets in a list $\mathcal{L}_s$ to be processed later. The details of the depth-first approach are presented in Algorithm \ref{DFCS}.

\begin{algo}
\begin{pseudocode}[shadowbox]{DepthFirstCanonicalSubsets}{s}
\label{DFCS}
\GLOBAL {G=\{g_1=e,g_2,\ldots,g_n\},\Aut(G)} \vspace{.1in} \\
\PROCEDURE{Extend}{S,k}
\COMMENT{$S \subseteq G$ and $|S| = k \leq s$}\\
\IF k = s \CTHEN  
\text{{proceed to step 2} ($S$ is a canonical $s$-subset of $G$)}
\ELSE 
\BEGIN
\text{denote the maximum element of $S$ by $g_m$}\\
\FOR i \GETS m+1 \TO n \DO \BEGIN
 S' \GETS S \cup \{g_i\}\\
 \IF \text{$S'$ is canonical} \CTHEN
   \mathsf{Extend}(S',k+1)\\
\END \END
\ENDPROCEDURE
\MAIN
\mathsf{Extend}(\{g_1\},1)\\
 \COMMENT{this is the initial call to $\mathsf{Extend}$}
\ENDMAIN
\end{pseudocode}
\end{algo}

The breadth-first method (Algorithm \ref{BFCS}) requires more storage than the depth-first method (Algorithm \ref{DFCS}). However, Algorithm \ref{BFCS} might be faster than {Algorithm \ref{DFCS}} because there is no recursion required in a breadth-first search.

\subsection{Finding all factorizations for a given group $G$}
\label{Sec3.4}

We now consider the problem of finding all factorizations of a given abelian group $G$. We use the following approach.
We first generate all the basic parameter sets. Then, for each basic parameter set, we  determine the inequivalent factorizations using the techniques presented in Sections \ref{canonicaltest.sec}--\ref{canon.sec}. However, this 
would be somewhat inefficient, as canonical $s$-subsets would be unnecessarily recomputed. We discuss in this section an enhanced algorithm that avoids this problem. 

First, generating all the basic parameter sets $(s,t,\lambda,\mu)$ for a given group $G$ with $|G| = n$ is straightforward, based on the conditions given in (\ref{eqn-basic}) and the remarks in the paragraph following (\ref{eqn-basic}). We have the following requirements:
\begin{itemize}
\item $2 \leq s \leq \lfloor \frac{n}{2} \rfloor$, 
\item $1 \leq \lambda , \mu \leq s-1$ and 
\item $\lambda \neq \mu$.
\end{itemize}
 Given $n, s, \lambda$ and $\mu$, the value of $t$ is
then determined by  equation (\ref{eqn-counting}) to be $t = (\lambda (n-1) + \mu)/s$. If $t$ is an integer such that $s \leq t \leq \lfloor \frac{n}{2} \rfloor$, then $(s,t,\lambda,\mu)$ is a basic parameter set. 

For $2 \leq j \leq \lfloor \frac{n}{2} \rfloor$, define 
\[ \mathcal{P}_j = \{ (s,t,\lambda,\mu)  : (s,t,\lambda,\mu)  \text{ is a basic parameter set with } s = j\} .\]
We present pseudocode to compute the sets $\mathcal{P}_j$ for $2 \leq j \leq \lfloor \frac{n}{2} \rfloor$ {in Algorithm \ref{BPS}}.

\begin{algo}
\begin{pseudocode}[shadowbox]{BasicParameterSets}{n}
\label{BPS}
\COMMENT{$|G| = n$}\\
\FOR s \GETS 2 \TO \lfloor \frac{n}{2} \rfloor \DO
\BEGIN
\mathcal{P}_s \GETS \emptyset\\
\FOR \lambda \GETS 1 \TO s-1 \DO
\BEGIN \FOR \mu \GETS 1 \TO s-1 \DO
\BEGIN 
\IF \lambda \neq \mu \THEN
\BEGIN
t \GETS \frac{\lambda (n-1) + \mu}{s}\\
\IF $t$ \text{ is an integer } \AND s \leq t \leq \lfloor \frac{n}{2} \rfloor \THEN
\mathcal{P}_s \GETS \mathcal{P}_s \cup \{ (s,t,\lambda,\mu ) \} 
\END 
\END
\END
\END\\
\RETURN{\mathcal{P}_2, \dots , \mathcal{P}_{\lfloor \frac{n}{2} \rfloor}}
\end{pseudocode}
\end{algo}

Having determined all the basic parameter sets, we then search for mates, starting with $\mathcal{P}_2$, then proceeding to $\mathcal{P}_3, \mathcal{P}_4$, etc. Observe that we need to generate all the canonical $j$-subsets
of $G$ in order to examine each of the basic parameter sets in $\mathcal{P}_j$. 

Let $s^*$ denote the largest value of $s$ such that $\mathcal{P}_s \neq \emptyset$ 
(of course $s^* \leq \lfloor \frac{n}{2} \rfloor$). 
Suppose we run Algorithm \ref{BFCS} with $s = s^*$. This algorithm returns all the lists $\mathcal{L}_1 , \dots , \mathcal{L}_{s^*}$. It is important to observe that we only run Algorithm \ref{BFCS} once. 

We let $j$ range from $2$ to $s^*$. For each value of $j$, we look at all the parameter sets in $\mathcal{P}_j$. Then, for each parameter set, we try to find mates for all the $j$-subsets in $\mathcal{L}_j$. 
Algorithm \ref{FIF} puts all the pieces together.

\begin{algo}
\begin{pseudocode}[shadowbox]{FindAllFactorizations}{G}
\label{FIF}
\COMMENT{$|G| = n$}\\
(\mathcal{P}_2, \dots , \mathcal{P}_{\lfloor \frac{n}{2} \rfloor}) \GETS
\mathsf{BasicParameterSets}(n)\\
s^* \GETS \max \{ s : \mathcal{P}_s \neq \emptyset \}\\
(\mathcal{L}_{1}, \dots , \mathcal{L}_{s^*}) \GETS \mathsf{BreadthFirstCanonicalSubsets}(s^*)\\
\FOR s \GETS 2 \TO {s^*} \DO
\BEGIN
\FORALL (s,t,\lambda,\mu)  \in {\mathcal{P}_{s}} \DO
\BEGIN
\FORALL S \in {\mathcal{L}_{s}} \DO
\BEGIN
\IF S \text{ has a mate } T \THEN
\BEGIN
(S',T') \GETS \mathsf{ComputeCanonicalFactorization}(S,T)\\
\IF (S',T') = (S,T) \CTHEN \OUTPUT {S,T}
\END
\END
\END
\END
\end{pseudocode}
\end{algo}

Algorithm \ref{FIF} is based on the breadth-first method of constructing canonical subsets (Algorithm \ref{BFCS}). As we have already mentioned, if storage is an issue, then it might be preferable to use the depth-first approach (Algorithm \ref{DFCS}). It is relatively straightforward to modify Algorithm \ref{FIF} to create a depth-first variation; we leave the details to the reader.

\section{Classifications of basic $(\la,\mu)$-fold factorizations of abelian groups of order at most $34$}
\label{sec-computation}

Applying the computational approach detailed in Section~\ref{sec-inequivalent}, up to equivalence, we are able to classify basic $(\la,\mu)$-fold factorizations of abelian groups of order at most $34$. We provide the examples in Table \ref{Exists} and nonexistence results in Table \ref{DNE}. 

All computations were performed on an Apple MacBook Pro equipped with an M1 chip and 16~GB of unified memory. The primary algorithms were implemented in C, compiled using Apple Clang (via the GCC frontend), and parallelized to execute concurrently across all 8 cores of the M1 processor. To further maximize computational efficiency, sets were represented as bit vectors using 64-bit unsigned integers  (\texttt{uint64\_t}), enabling the use of native bitwise operations for rapid set manipulation. To verify the accuracy and integrity of the C implementation, the computations were also independently developed in Maple~2024 and executed on the smaller example cyclic cases presented in this paper. Because Maple is roughly $100$ to $1,000$ times slower than C for these tasks, larger computations were performed exclusively in C.

\begin{center}
\renewcommand{\arraystretch}{1.2}
\begin{longtable}{@{}rlrrrr@{\hspace{12pt}}R{10.00cm}@{}}
\caption{Table of inequivalent basic $(\lambda,\mu)$-fold $(s,t)$-factorizations.}\label{Exists}\\ \hline
$n$ & group & $s$ & $t$ & $\lambda$  &$\mu$ & Authority  \\
\hline
\endfirsthead
\caption*{Table~\ref{Exists} (continued)}\\ \hline
$n$ & group & $s$ & $t$ & $\lambda$  &$\mu$ & Authority  \\
\hline 
\endhead
\endfoot
8&$\zed_{8} $&3&3&1&2&\cite[Lemma 4.1]{KLS-part1}\\ \hline
11&$\zed_{11} $&3&4&1&2&\cite[Lemma 4.1]{KLS-part1}\\ \hline
12&$\zed_{12} $&5&5&2&3&\cite[Theorem 4.2]{KLS-part1}\\ \hline
13&$\zed_{13} $&5&5&2&1&\cite[Theorem 4.4]{KLS-part1}\\ \hline
14&$\zed_{14} $&3&5&1&2&\cite[Lemma 4.1]{KLS-part1}\\ \hline
15&$\zed_{15} $&4&4&1&2&\cite[Lemma 4.1]{KLS-part1}\\ \hline
15&$\zed_{15} $&4&4&1&2&$S=\{0,1,4,5\}$,\allowbreak~$T=\{0,2,8,10\}$ \newline applying \cite[Theorem 4.5]{KLS-part1} to $(m_i)_{i=1}^4=(2,2,2,2)$ yields $(S,T)$ \\ \hline
16&$\zed_{16} $&7&7&3&4&\cite[Theorem 4.2]{KLS-part1}\\ \hline
17&$\zed_{17} $&3&6&1&2&\cite[Lemma 4.1]{KLS-part1}\\ \hline
17&$\zed_{17} $&5&7&2&3&\cite[Theorem 4.2]{KLS-part1}\\ \hline
18&$\zed_{18} $&5&7&2&1&\cite[Theorem 4.4]{KLS-part1}\\ \hline
19&$\zed_{19} $&4&5&1&2&\cite[Lemma 4.1]{KLS-part1}\\ \hline
19&$\zed_{19} $&7&8&3&2& $S=\{ 0,1,2,3,4,5,12\}$,\allowbreak~$T=\{ 2,3,6,8,11,12,16,17\}$\\ \hline
19&$\zed_{19} $&7&8&3&2&\cite[Theorem 4.4]{KLS-part1}\\ \hline
20&$\zed_{20} $&3&7&1&2&\cite[Lemma 4.1]{KLS-part1}\\ \hline
20&$\zed_{20} $&9&9&4&5&\cite[Theorem 4.2]{KLS-part1}\\ \hline
21&$\zed_{21} $&8&8&3&4&\cite[Theorem 4.2]{KLS-part1}\\ \hline
21&$\zed_{21} $&8&8&3&4&$S=\{0,1,2,3,5,6,7,8\}$,\allowbreak~$T=\{3,6,7,10,14,16,18,20\}$  \newline applying \cite[Theorem 4.5]{KLS-part1} to $(m_i)_{i=1}^4=(4,4,2,2)$ yields \newline $(S-5,T+5)$ \\ \hline
22&$\zed_{22} $&5&9&2&3&\cite[Theorem 4.2]{KLS-part1}\\ \hline
23&$\zed_{23} $&3&8&1&2&\cite[Lemma 4.1]{KLS-part1}\\ \hline
23&$\zed_{23} $&4&6&1&2&\cite[Lemma 4.1]{KLS-part1}\\ \hline
23&$\zed_{23} $&4&6&1&2&$S=\{0,1,4,5\}$,\allowbreak~$T=\{0,2,8,10,16,18\}$ \newline applying \cite[Theorem 4.5]{KLS-part1} to $(m_i)_{i=1}^4=(2,2,2,3)$ yields $(S,T)$ \\ \hline
23&$\zed_{23} $&5&9&2&1&\cite[Theorem 4.4]{KLS-part1}\\ \hline
23&$\zed_{23} $&7&10&3&4&\cite[Theorem 4.2]{KLS-part1}\\ \hline
24&$\zed_{24} $&5&5&1&2&\cite[Lemma 4.1]{KLS-part1}\\ \hline
24&$\zed_{24} $&7&7&2&3&\cite[Theorem 4.2]{KLS-part1}\\ \hline
24&$\zed_{24} $&11&11&5&6&\cite[Theorem 4.2]{KLS-part1}\\ \hline
25&$\zed_{25} $&7&7&2&1&\cite[Theorem 4.4]{KLS-part1}\\ \hline
25&$\zed_{25} $&9&11&4&3&\cite[Theorem 4.4]{KLS-part1}\\ \hline
25&$\zed_{25} $&9&11&4&3&$\begin{array}{@{}l@{}l@{}}S&=\{ 0,1,2,3,4,5,6,8,23\}\\T&=\{ 3,4,5,9,10,14,15,16,21,22,23\}\end{array}$\\ \hline
25&$\zed_{5}^2$&7&7&2&1&\cite[Example 4.20]{KLS-part1}\\ \hline
26&$\zed_{26} $&3&9&1&2&\cite[Lemma 4.1]{KLS-part1}\\ \hline
26&$\zed_{26} $&7&11&3&2&\cite[Theorem 4.4]{KLS-part1}\\ \hline
27&$\zed_{27} $&4&7&1&2&\cite[Lemma 4.1]{KLS-part1}\\ \hline
27&$\zed_{27} $&5&11&2&3&\cite[Theorem 4.2]{KLS-part1}\\ \hline
27&$\zed_{27} $&8&10&3&2&\cite[Theorem 4.4]{KLS-part1}\\ \hline
27&$\zed_{27} $&8&10&3&2&$S=\{ 0,1,2,3,7,8,9,10\}$, \allowbreak~$T=\{ 2,6,7,10,11,15,19,21,23,25 \}$ \newline applying \cite[Theorem 4.5]{KLS-part1} to $(m_i)_{i=1}^4=(4,5,2,2)$ yields \newline $(S-7,T+7)$ \\ \hline
27&$\zed_{27} $&8&10&3&2&$S=\{ 0,1,2,3,8,9,10,11 \}$, \allowbreak~$T=\{ 1,5,6,10,11,15,17,21,22,26 \}$ \newline applying \cite[Theorem 4.5]{KLS-part1} to $(m_i)_{i=1}^4=(4,2,2,5)$ yields $(S,T)$   \\ \hline
28&$\zed_{28} $&5&11&2&1&\cite[Theorem 4.4]{KLS-part1}\\ \hline
28&$\zed_{28} $&13&13&6&7&\cite[Theorem 4.2]{KLS-part1}\\ \hline
29&$\zed_{29} $&3&10&1&2&\cite[Lemma 4.1]{KLS-part1}\\ \hline
29&$\zed_{29} $&5&6&1&2&\cite[Lemma 4.1]{KLS-part1}\\ \hline
29&$\zed_{29} $&8&11&3&4&\cite[Theorem 4.2]{KLS-part1}\\ \hline
29&$\zed_{29} $&8&11&3&4&$\begin{array}{@{}l@{}l@{}}S&=\{ 0,1,2,3,11,12,13,14\}\\T&=\{ 0,3,7,8,12,15,16,20,22,24,28\}\end{array}$\\ \hline
29&$\zed_{29} $&8&11&3&4&$\begin{array}{@{}l@{}l@{}}S&=\{ 0,1,2,3,7,8,9,10\}\\T&=\{ 3,7,8,11,12,16,20,22,24,26,28\}\end{array}$\\ \hline
29&$\zed_{29} $&9&13&4&5&\cite[Theorem 4.2]{KLS-part1}\\ \hline
29&$\zed_{29} $&9&13&4&5&$\begin{array}{@{}l@{}l@{}}S&=\{ 0,1,2,3,4,5,6,7,18\}\\T&=\{ 3,4,6,9,11,13,16,18,19,23,25,26,28\}\end{array}$\\ \hline
29&$\zed_{29} $&9&13&4&5&$\begin{array}{@{}l@{}l@{}}S&=\{ 0,1,2,3,4,5,6,8,27\}\\T&=\{ 2,3,7,8,9,14,15,16,20,21,25,26,27\}\end{array}$\\ \hline
29&$\zed_{29} $&9&13&4&5&$\begin{array}{@{}l@{}l@{}}S&=\{ 0,1,2,5,6,7,17,18,19\}\\T&=\{ 0,1,4,7,8,11,14,15,18,21,22,24,27\}\end{array}$\\ \hline
29&$\zed_{29} $&12&12&5&4&\cite[Theorem 4.4]{KLS-part1}\\ \hline
29&$\zed_{29} $&12&12&5&4&$S=\{ 0,1,2,3,4,5,6,7,9,11,25,27 \}$ \newline $T=\{ 0,1,2,3,9,10,12,13,19,20,21,22 \}$ \newline applying \cite[Theorem 4.5]{KLS-part1} to $(m_i)_{i=1}^5=(2,6,2,2,3)$ yields $(T,S)$   \\ \hline
29&$\zed_{29} $&12&12&5&4&$S=\{ 0,1,2,3,4,5,7,9,11,23,25,27 \}$, \newline $T=\{ 2,3,4,5,7,8,16,17,19,20,21,22 \}$ \newline applying \cite[Theorem 4.5]{KLS-part1} to $(m_i)_{i=1}^5=(2,6,3,2,2)$ yields \newline $(T-7,S+7)$  \\ \hline
30&$\zed_{30} $&7&13&3&4&\cite[Theorem 4.2]{KLS-part1}\\ \hline
30&$\zed_{30} $&11&11&4&5&\cite[Theorem 4.2]{KLS-part1}\\ \hline
31&$\zed_{31} $&4&8&1&2&\cite[Lemma 4.1]{KLS-part1}\\ \hline
31&$\zed_{31} $&4&8&1&2&$S=\{0,1,4,5\}$,\allowbreak~$T=\{0,2,8,10,16,18,24,26\}$ \newline applying \cite[Theorem 4.5]{KLS-part1} to $(m_i)_{i=1}^4=(2,2,2,4)$ yields $(S,T)$ \\ \hline
31&$\zed_{31} $&7&9&2&3&\cite[Theorem 4.2]{KLS-part1}\\ \hline
31&$\zed_{31} $&11&14&5&4&\cite[Theorem 4.4]{KLS-part1}\\ \hline
31&$\zed_{31} $&11&14&5&4&$\begin{array}{@{}l@{}l@{}}S&=\{ 0,1,2,3,4,5,6,7,8,9,20\}\\T&=\{ 2,3,6,7,10,12,15,16,19,20,24,25,28,29\}\end{array}$\\ \hline
31&$\zed_{31} $&12&13&5&6&\cite[Theorem 4.2]{KLS-part1}\\ \hline
31&$\zed_{31} $&12&13&5&6&$\begin{array}{@{}l@{}l@{}}S&=\{ 0,1,2,3,4,5,8,9,10,11,12,13\}\\T&=\{ 0,2,3,6,9,12,15,16,18,21,22,27,28\}\end{array}$\\ \hline
31&$\zed_{31} $&12&13&5&6&$\begin{array}{@{}l@{}l@{}}S&=\{ 0,1,2,3,5,6,8,11,23,26,28,29\}\\T&=\{ 0,1,8,12,13,14,15,16,20,27,28,29,30\}\end{array}$\\ \hline
32&$\zed_{32} $&3&11&1&2&\cite[Lemma 4.1]{KLS-part1}~\\ \hline
32&$\zed_{32} $&5&13&2&3&\cite[Theorem 4.2]{KLS-part1}~\\ \hline
32&$\zed_{32} $&7&9&2&1&\cite[Theorem 4.4]{KLS-part1}~\\ \hline
32&$\zed_{32} $&15&15&7&8&\cite[Theorem 4.2]{KLS-part1}~\\ \hline
32&$\zed_{16} \times \zed_{2} $&14&14&6&10&$\begin{array}{@{}l@{}l@{}}S&=\left\{\begin{array}{@{}l@{}} (0,0),(1,0),(2,0),(3,0),(4,0),(5,0),(6,0),\\ (8,0),(1,1),(3,1),(5,1),(10,1),(12,1),(14,1)\end{array}\right\}\\[3ex]T&=\left\{\begin{array}{@{}l@{}} (0,0),(4,0),(5,0),(8,0),(10,0),(11,0),(14,0),\\ (15,0),(2,1),(5,1),(6,1),(11,1),(12,1),(15,1)\end{array}\right\}\end{array}$\\ \hline
32&$\zed_{16} \times \zed_{2} $&14&14&6&10&$\begin{array}{@{}l@{}l@{}}S&=\left\{\begin{array}{@{}l@{}} (0,0),(1,0),(2,0),(3,0),(4,0),(7,0),(9,0),\\ (13,0),(0,1),(2,1),(4,1),(5,1),(11,1),(15,1)\end{array}\right\}\\[3ex]T&=\left\{\begin{array}{@{}l@{}} (0,0),(6,0),(7,0),(9,0),(11,0),(12,0),(13,0),\\ (15,0),(0,1),(1,1),(3,1),(5,1),(6,1),(12,1)\end{array}\right\}\end{array}$\\ \hline
32&$\zed_{8} \times \zed_{4} $&9&11&3&6&$\begin{array}{@{}l@{}l@{}}S&=\left\{\begin{array}{@{}l@{}} (0,0),(1,0),(2,0),(0,1),\\ (1,1),(7,1),(2,2),(3,2),(4,2)\end{array}\right\}\\[3ex]T&=\left\{\begin{array}{@{}l@{}} (0,0),(3,0),(6,0),(0,1),(4,1),\\ (1,2),(4,2),(6,2),(1,3),(4,3),(7,3)\end{array}\right\}\end{array}$\\ \hline
32&$\zed_{8} \times \zed_{4} $&14&14&6&10&$\begin{array}{@{}l@{}l@{}}S&=\left\{\begin{array}{@{}l@{}} (0,0),(1,0),(2,0),(3,0),(4,0),(0,1),(1,1), \\ (3,1),(2,2),(5,2),(7,2),(1,3),(3,3),(4,3) \end{array}\right\}\\[3ex]T&=\left\{\begin{array}{@{}l@{}} (0,0),(2,0),(4,0),(5,0),(7,0),(0,1),(5,1), \\ (7,1),(1,2),(2,2),(3,2),(4,3),(5,3),(7,3) \end{array}\right\}\end{array}$ \newline applying \cite[Theorem 4.17]{KLS-part1} to $(S_1,T_1)=(\{0,3,6\},\{0,1,2\})$ in $\zed_{8}$ and $(S_2,T_2)=(\{0\},\{0\})$ in $\zed_{4}$ yields $(S_3,T_3)$, and $(S,T)=(\phi(S_3)+(5,3),\phi(T_3)-(5,3))$, where $\phi(x, y)=(3x+2y, x+3y)$ for $(x,y) \in \zed_8 \times \zed_4$\\ \hline
32&$\zed_{8} \times \zed_{4} $&14&14&6&10&$\begin{array}{@{}l@{}l@{}}S&=\left\{\begin{array}{@{}l@{}} (0,0),(1,0),(2,0),(3,0),(0,1),(1,1),(2,1),\\ (5,1),(0,2),(1,2),(2,2),(7,2),(4,3),(6,3)\end{array}\right\}\\[3ex]T&=\left\{\begin{array}{@{}l@{}} (0,0),(1,0),(3,0),(6,0),(2,1),(4,1),(0,2),\\ (3,2),(5,2),(6,2),(0,3),(3,3),(6,3),(7,3)\end{array}\right\}\end{array}$\\ \hline
32&$\zed_{8} \times \zed_{4} $&14&14&6&10&$\begin{array}{@{}l@{}l@{}}S&=\left\{\begin{array}{@{}l@{}} (0,0),(1,0),(2,0),(3,0),(0,1),(1,1),(3,1),\\ (4,1),(1,2),(2,2),(3,2),(4,2),(5,3),(7,3)\end{array}\right\}\\[3ex]T&=\left\{\begin{array}{@{}l@{}} (2,0),(4,0),(5,0),(7,0),(1,1),(3,1),(0,2),\\ (2,2),(5,2),(7,2),(0,3),(4,3),(5,3),(7,3)\end{array}\right\}\end{array}$\\ \hline
32&$\zed_{8} \times \zed_{4} $&14&14&6&10&$\begin{array}{@{}l@{}l@{}}S&=\left\{\begin{array}{@{}l@{}} (0,0),(1,0),(2,0),(3,0),(4,0),(0,1),(1,1),\\ (2,1),(1,2),(2,2),(3,2),(4,3),(5,3),(6,3)\end{array}\right\}\\[3ex]T&=\left\{\begin{array}{@{}l@{}} (0,0),(2,0),(4,0),(5,0),(7,0),(2,1),(4,1),\\ (7,1),(2,2),(5,2),(7,2),(0,3),(3,3),(6,3)\end{array}\right\}\end{array}$\\ \hline
32&$\zed_{8} \times \zed_{4} $&14&14&6&10&$\begin{array}{@{}l@{}l@{}}S&=\left\{\begin{array}{@{}l@{}} (0,0),(1,0),(2,0),(3,0),(4,0),(0,1),(1,1),\\ (7,1),(5,2),(6,2),(7,2),(0,3),(1,3),(7,3)\end{array}\right\}\\[3ex]T&=\left\{\begin{array}{@{}l@{}} (0,0),(2,0),(4,0),(5,0),(7,0),(1,1),(4,1),\\ (7,1),(1,2),(3,2),(6,2),(1,3),(4,3),(7,3)\end{array}\right\}\end{array}$\\ \hline
32&$\zed_{8} \times \zed_{2}^{2} $&14&14&6&10&$\begin{array}{@{}l@{}l@{}}S&=\left\{\begin{array}{@{}l@{}} (0,0,0),(1,0,0),(2,0,0),(3,0,0),(4,0,0),(1,1,0),(2,1,0), \\ (3,1,0),(1,0,1),(2,0,1),(3,0,1),(5,1,1),(6,1,1),(7,1,1) \end{array}\right\}\\[3ex]T&=\left\{\begin{array}{@{}l@{}} (0,0,0),(2,0,0),(4,0,0),(5,0,0),(7,0,0),(2,1,0),(5,1,0), \\ (7,1,0),(2,0,1),(5,0,1),(7,0,1),(1,1,1),(3,1,1),(6,1,1) \end{array}\right\}\end{array}$ \newline applying \cite[Theorem 4.17]{KLS-part1} to $(S_1,T_1)=(\{0,3,6\},\{0,1,2\})$ in $\zed_{8}$ and $(S_2,T_2)=(\{(0,0)\},\{(0,0)\})$  in $\zed_{2}^2$ yields $(S_3,T_3)$, and $(S,T)=(\phi(S_3)+(5,0,0),\phi(T_3)-(5,0,0))$, where $\phi(x,y,z) = (3x + 4y + 4z,y,z)$ for $(x,y,z) \in \zed_8 \times \zed_2^2$\\ \hline 
32&$\zed_{8} \times \zed_{2}^{2} $&14&14&6&10&$\begin{array}{@{}l@{}l@{}}S&=\left\{\begin{array}{@{}l@{}} (0,0,0),(1,0,0),(2,0,0),(3,0,0),(0,1,0),(1,1,0),(2,1,0),\\ (7,1,0),(0,0,1),(1,0,1),(5,0,1),(6,0,1),(2,1,1),(4,1,1)\end{array}\right\}\\[3ex]T&=\left\{\begin{array}{@{}l@{}} (0,0,0),(1,0,0),(3,0,0),(6,0,0),(0,1,0),(3,1,0),(5,1,0),\\ (6,1,0),(0,0,1),(2,0,1),(3,0,1),(7,0,1),(4,1,1),(6,1,1)\end{array}\right\}\end{array}$\\ \hline
32&$\zed_{8} \times \zed_{2}^{2} $&14&14&6&10&$\begin{array}{@{}l@{}l@{}}S&=\left\{\begin{array}{@{}l@{}} (0,0,0),(1,0,0),(2,0,0),(3,0,0),(4,0,0),(0,1,0),(1,1,0),\\ (2,1,0),(0,0,1),(1,0,1),(2,0,1),(5,1,1),(6,1,1),(7,1,1)\end{array}\right\}\\[3ex]T&=\left\{\begin{array}{@{}l@{}} (0,0,0),(2,0,0),(4,0,0),(5,0,0),(7,0,0),(0,1,0),(3,1,0),\\ (6,1,0),(0,0,1),(3,0,1),(6,0,1),(1,1,1),(3,1,1),(6,1,1)\end{array}\right\}\end{array}$\\ \hline
32&$\zed_{8} \times \zed_{2}^{2} $&14&14&6&10&$\begin{array}{@{}l@{}l@{}}S&=\left\{\begin{array}{@{}l@{}} (0,0,0),(1,0,0),(2,0,0),(3,0,0),(4,0,0),(0,1,0),(1,1,0),\\ (7,1,0),(0,0,1),(1,0,1),(7,0,1),(5,1,1),(6,1,1),(7,1,1)\end{array}\right\}\\[3ex]T&=\left\{\begin{array}{@{}l@{}} (0,0,0),(2,0,0),(4,0,0),(5,0,0),(7,0,0),(1,1,0),(4,1,0),\\ (7,1,0),(1,0,1),(4,0,1),(7,0,1),(1,1,1),(3,1,1),(6,1,1)\end{array}\right\}\end{array}$\\ \hline
33&$\zed_{33}$ & 5 & 13 & 2 & 1 &\cite[Theorem 4.4]{KLS-part1}~\\ \hline
33&$\zed_{33}$ & 7 & 14 & 3 & 2 &\cite[Theorem 4.4]{KLS-part1}~\\ \hline
33&$\zed_{33}$ & 10 & 10 & 3 & 4 &\cite[Theorem 4.2]{KLS-part1}~\\ \hline
33&$\zed_{33}$ & 10 & 10 & 3 & 4 & $S=\{ 0,1,2,3,4,8,9,10,11,12 \}$ \newline $T=\{ 2,7,9,12,14,19,24,25,29,30 \}$ \newline applying \cite[Theorem 4.5]{KLS-part1} to $(m_i)_{i=1}^4=(5,5,2,2)$ yields \newline $(S-8,T+8)$  \\ \hline
34&$\zed_{34}$ & 5 & 7 & 1 & 2 &\cite[Lemma 4.1]{KLS-part1}~\\ \hline
34&$\zed_{34}$ & 9 & 15 & 4 & 3 &\cite[Theorem 4.4]{KLS-part1}~\\ \hline
34&$\zed_{34}$ & 13 & 13 & 5 & 4 &\cite[Theorem 4.4]{KLS-part1}~\\ \hline
34&$\zed_{34}$ & 13 & 13 & 5 & 4 &$\begin{array}{@{}l@{}l@{}}S&=\{ 0,1,2,3,4,6,12,15,17,21,23,26,32 \}\\T&=\{ 2,6,7,9,14,15,16,21,23,24,28,31,33 \}\end{array}$\\ \hline
\end{longtable}
\end{center}

{\RaggedRight
\setlength{\LTcapwidth}{\textwidth}
\renewcommand{\arraystretch}{1.15}

\begin{ThreePartTable}

\begin{TableNotes}
\normalsize
\item[]
\item[]{\leftskip=-\leftmargin Table~\ref{DNE} notes.\vspace{0.1in}\par}
\item[]{\leftskip=-\leftmargin A superscript records that nonexistence already follows from \cite[Section 5]{KLS-part1}:\vspace{0.1in}\par}
\item[a] \cite[Theorem 5.1(3)]{KLS-part1},
\item[b] \cite[Theorem 5.2(3)]{KLS-part1},
\item[c] \cite[Corollary 5.3]{KLS-part1}.\vspace{0.1in}
\item[]{\leftskip=-\leftmargin When more than one result applies, the most specific one is cited. Parameter sets carrying no superscript are those whose nonexistence follows only from the exhaustive computer search.\par}
\end{TableNotes}

\begin{longtable}{@{}r@{\hspace{4pt}}|@{\hspace{4pt}}p{0.19\textwidth}@{\hspace{4pt}}|@{\hspace{4pt}}p{0.715\textwidth}@{}}
\caption{Basic parameter sets $(s,t,\la,\mu)$ for which no $(\la,\mu)$-fold factorization of the abelian group $G$ of order $n \le 34$ exists. }\label{DNE}\\
\hline
$n$&group $G$&$(s,t,\la,\mu)$\\
\hline
\endfirsthead
\caption*{Table \ref{DNE} (continued)}\\ \hline
$n$&group $G$&$(s,t,\la,\mu)$\\
\hline
\endhead
\hline
\endfoot
\hline
    \insertTableNotes
\endlastfoot
$8$
& $\zed_{4} \times \zed_{2}$
& $(3,3,1,2)$.\\

$8$
& $\zed_{2}^{3}$
& $(3,3,1,2)^{c}$.\\

$12$
& $\zed_{6} \times \zed_{2}$
& $(5,5,2,3)$.\\

$14$
& $\zed_{14}$
& $(4,4,1,3)^{a}$,
$(5,6,2,4)$.\\

$16$
& $\zed_{4}^2$,
 $\zed_{8} \times \zed_{2}$,
 $\zed_{4} \times \zed_{2}^{2}$
& $(7,7,3,4)$.\\

$16$
& $\zed_{2}^{4}$
& $(7,7,3,4)^{c}$.\\
$17$
& $\zed_{17}$
& $(6,6,2,4)$,
$(7,7,3,1)$.\\

$18$
& $\zed_{18}$
& $(4,5,1,3)$,
$(7,8,3,5)$.\\

$18$
& $\zed_{6} \times \zed_{3}$
& $(4,5,1,3)$,
$(5,7,2,1)$,
$(7,8,3,5)$.\\

$19$
& $\zed_{19}$
& $(5,8,2,4)$.\\

$20$
& $\zed_{10} \times \zed_{2}$
& $(3,7,1,2)$,
$(6,7,2,4)$,
$(7,9,3,6)$,
$(8,8,3,7)$,
$(9,9,4,5)$.\\

$20$
& $\zed_{20}$
& $(6,7,2,4)$,
$(7,9,3,6)$,
$(8,8,3,7)$.\\

$21$
& $\zed_{21}$
& $(9,9,4,1)$.\\

$22$
& $\zed_{22}$
& $(4,6,1,3)^{a}$,
$(5,5,1,4)$,
$(8,8,3,1)^{a}$,
$(9,10,4,6)$.\\

$23$
& $\zed_{23}$
& $(5,5,1,3)$,
$(6,8,2,4)$,
$(7,7,2,5)$,
$(8,9,3,6)$,
$(9,10,4,2)$.\\

$24$
& $\zed_{12} \times \zed_{2}$
& $(5,5,1,2)$,
$(5,10,2,4)$,
$(7,7,2,3)$,
$(7,10,3,1)$,
$(9,11,4,7)$,
$(10,10,4,8)$,
$(11,11,5,6)$.\\

$24$
& $\zed_{6} \times \zed_{2}^{2}$
& $(5,5,1,2)^{b}$,
$(5,10,2,4)$,
$(7,7,2,3)$,
$(7,10,3,1)$,
$(9,11,4,7)$,
$(10,10,4,8)$,
$(11,11,5,6)^{b}$.\\

$24$
& $\zed_{24}$
& $(5,10,2,4)$,
$(7,10,3,1)$,
$(9,11,4,7)$,
$(10,10,4,8)$.\\

$25$
& $\zed_{25}$
& $(7,11,3,5)$,
$(11,11,5,1)$.\\

$25$
& $\zed_{5}^2$
& $(7,11,3,5)$,
$(9,11,4,3)$,
$(11,11,5,1)$.\\

$26$
& $\zed_{26}$
& $(4,7,1,3)$,
$(6,9,2,4)$,
$(7,8,2,6)$,
$(8,10,3,5)^{a}$,
$(9,9,3,6)$,
$(9,12,4,8)$,
$(11,12,5,7)$.\\

$27$
& $\zed_{27}$
& $(5,6,1,4)$,
$(7,8,2,4)$,
$(7,12,3,6)$,
$(10,11,4,6)$,
$(11,12,5,2)$.\\

$27$
& $\zed_{3}^3$,
 $\zed_{9} \times \zed_{3}$
& $(4,7,1,2)$,
$(5,6,1,4)$,
$(5,11,2,3)$,
$(7,8,2,4)$,
$(7,12,3,6)$,
$(8,10,3,2)$,
$(10,11,4,6)$,
$(11,12,5,2)$.\\

$28$
& $\zed_{14} \times \zed_{2}$
& $(5,6,1,3)$,
$(5,11,2,1)$,
$(8,11,3,7)$,
$(10,11,4,2)$,
$(11,13,5,8)$,
$(12,12,5,9)$,
$(13,13,6,7)$.\\

$28$
& $\zed_{28}$
& $(5,6,1,3)$,
$(8,11,3,7)$,
$(10,11,4,2)$,
$(11,13,5,8)$,
$(12,12,5,9)$.\\

$29$
& $\zed_{29}$
& $(5,12,2,4)$,
$(6,10,2,4)$,
$(9,10,3,6)$,
$(10,12,4,8)$,
$(11,11,4,9)$,
$(11,13,5,3)$,
$(13,13,6,1)$.\\

$30$
& $\zed_{30}$
& $(4,8,1,3)^{a}$,
$(7,9,2,5)$,
$(8,8,2,6)$,
$(8,11,3,1)$,
$(9,13,4,1)$,
$(11,14,5,9)$,
$(12,13,5,11)$,
$(13,14,6,8)$.\\

$31$
& $\zed_{31}$
& $(7,13,3,1)$,
$(8,8,2,4)$,
$(8,12,3,6)$,
$(9,14,4,6)$,
$(11,11,4,1)$,
$(13,14,6,2)$.\\

$32$
& $\zed_{16} \times \zed_{2}$
& $(3,11,1,2)$,
$(5,7,1,4)$,
$(5,13,2,3)$,
$(6,6,1,5)$,
$(6,11,2,4)$,
$(7,9,2,1)$,
$(7,14,3,5)$,
$(9,11,3,6)$,
$(9,14,4,2)$,
$(10,10,3,7)$,
$(10,13,4,6)$,
$(11,12,4,8)$,
$(11,15,5,10)$,
$(12,13,5,1)$,
$(13,15,6,9)$,
$(15,15,7,8)$.\\

$32$
& $\zed_{8} \times \zed_{4}$
& $(3,11,1,2)$,
$(5,7,1,4)$,
$(5,13,2,3)$,
$(6,6,1,5)$,
$(6,11,2,4)$,
$(7,9,2,1)$,
$(7,14,3,5)$,
$(9,14,4,2)$,
$(10,10,3,7)$,
$(10,13,4,6)$,
$(11,12,4,8)$,
$(11,15,5,10)$,
$(12,13,5,1)$,
$(13,15,6,9)$,
$(15,15,7,8)$.\\

$32$
& $\zed_{8} \times \zed_{2}^{2}$
& $(3,11,1,2)^{b}$,
$(5,7,1,4)$,
$(5,13,2,3)^{b}$,
$(6,6,1,5)$,
$(6,11,2,4)$,
$(7,9,2,1)$,
$(7,14,3,5)$,
$(9,11,3,6)$,
$(9,14,4,2)$,
$(10,10,3,7)$,
$(10,13,4,6)$,
$(11,12,4,8)$,
$(11,15,5,10)$,
$(12,13,5,1)$,
$(13,15,6,9)$,
$(15,15,7,8)$.\\

$32$
& $\zed_{4}^{2} \times \zed_{2}$
& $(3,11,1,2)^{b}$,
$(5,7,1,4)$,
$(5,13,2,3)^{b}$,
$(6,6,1,5)$,
$(6,11,2,4)$,
$(7,9,2,1)$,
$(7,14,3,5)$,
$(9,11,3,6)$,
$(9,14,4,2)$,
$(10,10,3,7)$,
$(10,13,4,6)$,
$(11,12,4,8)$,
$(11,15,5,10)$,
$(12,13,5,1)$,
$(13,15,6,9)$,
$(14,14,6,10)$,
$(15,15,7,8)$.\\

$32$
& $\zed_{4} \times \zed_{2}^{3}$
& $(3,11,1,2)^{b}$,
$(5,7,1,4)$,
$(5,13,2,3)^{b}$,
$(6,6,1,5)$,
$(6,11,2,4)$,
$(7,9,2,1)^{b}$,
$(7,14,3,5)^{b}$,
$(9,11,3,6)$,
$(9,14,4,2)^{b}$,
$(10,10,3,7)$,
$(10,13,4,6)$,
$(11,12,4,8)$,
$(11,15,5,10)$,
$(12,13,5,1)$,
$(13,15,6,9)$,
$(14,14,6,10)$,
$(15,15,7,8)$.\\

$32$
& $\zed_{32}$
& $(5,7,1,4)$,
$(6,6,1,5)$,
$(6,11,2,4)$,
$(7,14,3,5)$,
$(9,11,3,6)$,
$(9,14,4,2)$,
$(10,10,3,7)$,
$(10,13,4,6)$,
$(11,12,4,8)$,
$(11,15,5,10)$,
$(12,13,5,1)$,
$(13,15,6,9)$,
$(14,14,6,10)$.\\

$32$
& $\zed_{2}^{5}$
& $(3,11,1,2)^{c}$,
$(5,7,1,4)$,
$(5,13,2,3)^{c}$,
$(6,6,1,5)$,
$(6,11,2,4)^{c}$,
$(7,9,2,1)^{c}$,
$(7,14,3,5)^{c}$,
$(9,11,3,6)$,
$(9,14,4,2)^{c}$,
$(10,10,3,7)$,
$(10,13,4,6)^{c}$,
$(11,12,4,8)$,
$(11,15,5,10)$,
$(12,13,5,1)$,
$(13,15,6,9)$,
$(14,14,6,10)$,
$(15,15,7,8)^{c}$.\\

33& $\zed_{33}$
&$( 5, 7, 1, 3)$,
$( 6, 6, 1, 4)$,
$( 7,10, 2, 6)$,
$( 9,15, 4, 7)$,
$(10,13, 4, 2)$,
$(12,14, 5, 8)$,
$(13,13, 5, 9)$,
$(13,15, 6, 3)$,
$(14,14, 6, 4)$,
$(15,15, 7, 1)$.\\

34& $\zed_{34}$
&$( 4, 9, 1, 3)$,
$( 5,14, 2, 4)$,
$( 6, 6, 1, 3)^{a}$,
$( 7,10, 2, 4)$,
$( 7,15, 3, 6)$,
$( 8, 9, 2, 6)$,
$( 8,13, 3, 5)$,
$(10,10, 3, 1)^{a}$,
$(10,14, 4, 8)$,
$(12,14, 5, 3)^{a}$,
$(13,16, 6,10)$,
$(14,15, 6,12)$,
$(15,16, 7, 9)$.
\end{longtable}
\end{ThreePartTable}
}

\section*{AI usage disclosure}
 All underlying algorithms and mathematical methods were developed entirely by the authors. The primary computations reported in Tables \ref{Exists} and \ref{DNE} were performed using a C implementation, with Claude used strictly as an auxiliary tool to cross-verify the resulting data.


\begin{thebibliography}{0}

\bibitem{BHS}
G.~Bacs\'{o}, L.~H\'{e}thelyi and P.~Sziklai.  
New near-factorizations of finite groups.
\emph{Studia Sci. Math. Hungar.} {\bf 45} (2008), 493--510.

\bibitem{DB56}  
N.~G. de Bruijn.
On number systems.
\emph{Nieuw Arch. Wisk. (3)} {\bf 4} (1956), 15--17.

\bibitem{CGHK}
D.~de Caen, D.~A. Gregory, I.~G. Hughes and D.~L. Kreher.
Near-factors of finite groups.
\emph{Ars Combin.} {\bf 29} (1990), 53--63.

\bibitem{dummitfoote}
D.~S. Dummit and R.~M. Foote. \emph{Abstract Algebra}. 3rd ed. Hoboken, NJ: John Wiley \& Sons, 2004.

\bibitem{Hillar-Rhea}
C.~J. Hillar and D.~L. Rhea.
\newblock Automorphisms of finite abelian groups.
\newblock \emph{Amer. Math. Monthly} \textbf{114} (2007), 917--923.

\bibitem{KP}
P.~Kaski and O.~Pottonen.
\newblock {\em libexact User's Guide, Version 1.0}.
\newblock HIIT Technical Reports 2008-1, Helsinki Institute for Information
Technology HIIT, Helsinki, 2008.
\newblock Available at \url{https://pottonen.kapsi.fi/libexact.html}.

\bibitem{Knuth}
D.~E. Knuth.
\newblock Dancing links.
\newblock In J.~Davies, B.~Roscoe, and J.~Woodcock, editors,
{\em Millennial Perspectives in Computer Science}, pages 187--214.
\newblock Palgrave, Basingstoke, England, 2000.

\bibitem{KLS}
D.~L. Kreher, S.~Li, and D.~R. Stinson.
$\la$-fold near-factorizations of groups.
\newblock {\em Des. Codes Cryptogr.}, {\bf 94} (2026), Paper No. 106, 34 pp.

\bibitem{KLS-part1}
D.~L. Kreher, S.~Li, and D.~R. Stinson.
The theory of $(\la,\mu)$-fold factorizations of groups.
\newblock {\em arXiv preprint arXiv:2610.15826v2}, 2026.

\bibitem{KMS}
D.~L. Kreher, W.~J. Martin and D.~R. Stinson.
Uniqueness and explicit computation of mates in near-factorizations.
\newblock {\em J. Algebraic Combin.}, {\bf 63} (2026), Paper No. 53, 21 pp.

\bibitem{KPS2}
D.~L. Kreher, M.~B. Paterson and D.~R. Stinson.
Near-factorizations of dihedral groups.
\newblock {\em Des. Codes Cryptogr.}, {\bf 93} (2025), 5015--5038.

\bibitem{LM}
S.~Li and K.~Momihara.
Cyclotomic construction of $\la$-fold near-factorizations of cyclic groups.
\newblock {\em arXiv preprint arXiv:2507.18045}, 2025.


\bibitem{Mader}
A.~Mader.
\newblock The automorphism group of finite abelian {$p$}-groups.
\newblock {\em Ann. Sci. Math. Qu\'{e}bec}, {\bf 36} (2012), 559--575.

\bibitem{Pace}
N.~Pace.
One-overlapped factorizations of non-abelian groups.
\emph{Graphs Combin.} {\bf 34} (2018), 1459--1467.


\bibitem{Pech03}
A.~P\^{e}cher. 
Partitionable graphs arising from near-factorizations of finite groups.
{\em Discrete Math.} {\bf 269} (2003), 191--218.

\bibitem{Pech04}
A.~P\^{e}cher. 
Cayley partitionable graphs and near-factorizations of finite groups.
{\em Discrete Math.} {\bf 276} (2004), 295--311.

\bibitem{Read}
R.~C. Read.
Every one a winner or how to avoid isomorphism search when cataloguing combinatorial configurations.
\emph{Ann. Discrete Math.} {\bf 2} (1978), 107--120.

\bibitem{SS13}
T.~Sakuma and H.~Shinohara.
Krasner near-factorizations and 1-overlapped factorizations.
In \emph{The Seventh European Conference on Combinatorics, Graph Theory and
Applications}, {\it EuroComb 2013}, Pisa, 2013, pp.\ 391--395.

\bibitem{SS14}
T.~Sakuma and H.~Shinohara.
On circulant thin {L}ehman matrices.
\emph{J. Algebraic Combin.} {\bf 40} (2014), 939--959.

\bibitem{Shi}
H.~Shinohara.
Thin {L}ehman matrices arising from finite groups.
\emph{Linear Algebra Appl.} {\bf 436} (2012), 850--857.

\bibitem{Szabo04}
S.~Szab\'{o}.
\newblock {\em Topics in factorization of abelian groups}.
\newblock Birkh\"{a}user Verlag, Basel, 2004.

\bibitem{SS}
S.~Szab\'{o} and A.~D. Sands.
\newblock {\em Factoring groups into subsets}, volume 257 of
{\em Lecture Notes in Pure and Applied Mathematics}.
\newblock CRC Press, Boca Raton, FL, 2009.

\end{thebibliography}
\end{document}